\documentclass[12pt]{amsart}

\usepackage{amssymb}

\usepackage{anysize}
\marginsize{2cm}{2cm}{2cm}{2cm}

\usepackage{wrapfig }
\usepackage{graphicx}
\usepackage{caption}
\usepackage{subcaption}

\usepackage{versions}
\usepackage{hyperref}
\hypersetup{
  colorlinks   = true, 
  urlcolor     = blue, 
  linkcolor    = blue, 
  citecolor   = red 
}
\usepackage{cancel}
\usepackage{yhmath}
\usepackage{enumitem}
\usepackage{stackengine}

\newtheorem{theorem}{Theorem}[section]
\newtheorem*{theorem*}{Theorem}

\newtheorem*{proposition*}{Proposition}
\newtheorem{lemma}[theorem]{Lemma}
\newtheorem*{lemma*}{Lemma}

\newtheorem*{conjecture*}{Conjecture}
\newtheorem{claim}[theorem]{Claim}

\newtheorem*{corollary*}{Corollary}
\newtheorem*{observation*}{Observation}
\newtheorem*{problem*}{Problem}

\theoremstyle{definition}

\newtheorem*{definition*}{Definition}

\theoremstyle{remark}
\newtheorem{remark}[theorem]{Remark}
\newtheorem*{question*}{Question}

\numberwithin{equation}{section}

\DeclareMathOperator{\conv}{conv}
\DeclareMathOperator{\aff}{aff}

\newcommand*{\where}{\ \ifnum\currentgrouptype=16 \middle\fi|\ }

\renewcommand{\epsilon}{\varepsilon}
\renewcommand{\phi}{\varphi}
\renewcommand{\kappa}{\varkappa}
\renewcommand{\theta}{\vartheta}

\def\R{{\mathbb R}}

\def\F{{\mathbb F}}

\title{The affine Tverberg theorem revisited}

\author{Sergey Avvakumov}

\author{Roman Karasev}

\address{Sergey~Avvakumov,  School of Mathematical Sciences, Tel Aviv University, Tel Aviv 69978, Israel}
\email{savvakumov@gmail.com}

\address{Roman~Karasev, Institute for Information Transmission Problems RAS, Bolshoy Karetny per. 19, Moscow, Russia 127994}
\email{r\_n\_karasev@mail.ru}
\urladdr{http://www.rkarasev.ru/en/}


\begin{document}

\begin{abstract}
We prove the B\'ar\'any--Kalai--Tverberg conjecture extending the affine version of Tverberg's theorem to simplicial balls and polytopes other than the simplex:

For any affine map $\phi : P\to \mathbb R^d$ from a convex polytope $P$ of dimension $N=(d+1)(r-1)$ with $r\geq 2$,  there exist $r$ pairwise disjoint faces $F_1,\ldots,F_r\subset\partial P$ such that $\phi(F_1)\cap\dots\cap\phi(F_r)\neq\varnothing$. A similar statement holds for a simplicial ball $P$ with $\phi$ assumed affine on all of its faces.
\end{abstract}

\maketitle

\section{Introduction}

Tverberg's theorem \cite{tverberg1966} asserts that \emph{for any affine map $\phi : \Delta^N\to \mathbb R^d$ from the simplex of dimension $N=(d+1)(r-1)$ there exist $r$ pairwise disjoint faces $F_1,\ldots,F_r\subset\Delta^N$ such that $\phi(F_1)\cap\dots\cap\phi(F_r)\neq\varnothing$.}

In \cite{bss1981,ozaydin1987,volovikov1996} this theorem was extended to continuous maps $\phi$ under the assumption that $r$ is a prime power using topological tools of equivariant maps. The results in \cite{ozaydin1987,mabillard2015,gromov2010} implied, as was first noted by Florian Frick \cite[Section 2.9]{frick2015}, that the theorem for continuous $\phi$ may fail for $r$ not a prime power (see the survey \cite{skopenkov2018} for an exposition of these counterexamples).

Helge~Tverberg \cite{gruber1979problems}, and Imre~B\'ar\'any and Gil~Kalai conjectured \cite[Conjecture 4.5]{barany2022helly} that the theorem should still hold for spaces more general than the simplex, in particular, with $\Delta^N$ replaced by an arbitrary convex polytope. Counterexamples above apply in that case as well, so with continuous $\phi$ one has to assume $r$ to be a prime power. Thus, a version for simplicial $(N-1)$-spheres (instead of $\Delta^N$), with $\phi$ continuous, and $r$ a prime power was established in \cite{hktt2023}.

The main difficulty of the B\'ar\'any--Kalai--Tverberg conjecture, however, is that the standard equivariant methods are not applicable in the non-prime power case. In \cite{soberon-zerbib2026} it is claimed that the polytopal case of Theorem~\ref{theorem:main} immediately follows from Gr\"unbaum's classical theorem that any polytope is a refinement of a simplex; below in Remark~\ref{remark:grunbaum} we explain why we think that there is no such immediate implication.

\begin{theorem}[Tverberg's theorem for polytopes and simplicial balls]
\label{theorem:main}
Let $d\geq 1$ and $r\geq 2$ be integers such that
\[
N=(d+1)(r-1).
\]
Let $P$ be either an $N$-dimensional convex polytope (polytopal case) or an $N$-dimensional simplicial ball (simplicial-ball case).

Let
\[
\phi:P\longrightarrow\R^d
\]
be a continuous map which is affine on $P$ in the polytopal case and affine on every face of $P$ in the simplicial-ball case.

Then there are $r$ pairwise disjoint faces
$F_1,\ldots,F_r\subset \partial P$ such that
\[
\phi(F_1)\cap\cdots\cap\phi(F_r)\neq\varnothing.
\]
\end{theorem}

\begin{remark}
Comparing with the results in \cite{hktt2023}, we may start with a simplicial $(N-1)$-dimensional sphere $K$ and let $P$ be the cone over $K$. Then any map $\phi' : K\to \mathbb R^d$ affine on faces extends to $P$ and after application of Theorem~\ref{theorem:main} produces the Tverberg-type theorem for $K$.
\end{remark}

Our proof  uses a combination of topological and geometric methods. We start with Sarkaria's standard colorful Carath\'erodory approach \cite{sarkaria1992tverberg}. The two crucial facts we need to establish, the second of which is much easier in the simplicial case $P=\Delta^N$, are the following: the $r$-fold deleted join $D$ has vanishing $\widetilde H^{N-1}$ (Lemma~\ref{lemma:diagonal_acyclicity}); and moreover, if points $y\in S^{N-1}$ parametrize linear functionals on $Q$ which are zero on the diagonal, $X\subset D\times S^{N-1}$ is a set of pairs $(z,y)$ with the functional $y$ non-negative on the minimal face containing $z$, then there is a suitably defined subset $Y\subset X$ whose projection to $S^{N-1}$ has $(N-2)$-acyclic fibers, Lemma~\ref{lemma:Y}.

\begin{remark}
\label{remark:grunbaum}
In \cite{soberon-zerbib2026} it is claimed that the polytopal case of Theorem~\ref{theorem:main} immediately follows from the classical Tverberg theorem and Gr\"unbaum's theorem that any polytope is a refinement of a simplex. Gr\"unbaum's theorem specifically says \cite[Theorem 11.1.1, p.~200]{grunbaum2003} that there is a homeomorphism $\psi:P\to \Delta^N$ such that the preimage of any face of $\Delta^N$ is a subcomplex of $P$. A map $\phi:P\to\R^d$ with no $r$-tuple points then yields the map $\phi\circ\psi^{-1}:\Delta^N\to \R^d$ again with no $r$-tuple points. This map is not guaranteed to be affine though, even if $\phi$ is, and hence the classical Tverberg's theorem does not apply.
\end{remark}

\subsection*{AI use statement} We thank ChatGPT 5.6 for its help with the proof. More precisely, we credit ChatGPT with the construction of $Y$ in Section~\ref{section:Ry} and the statement and proof of Lemma~\ref{lemma:acyclic-gluing}. The authors independently verified all the arguments of the model. Most of the text is written by the authors; the proof of Lemma~\ref{lemma:acyclic-gluing} was initially written by the model and later revised by the authors. 

\subsection*{Lean formalization} The paper's proof of Theorem~\ref{theorem:main} was formalized and checked in Lean 4, see \url{https://github.com/savvakumov/AffineTverberg} for details. Like the proof, the formalization assumes $r\geq 3$ in the the simplicial-ball case, since $r=2$ (as well as $r$ equal any prime power) is already covered by \cite{hktt2023}.

\subsection*{Acknowledgments} We thank Gil Kalai and Uli Wagner for useful conversations about the problem.

\section{Standard preparations}

We follow the standard approach of passing from the map $\phi$ from the simplex (in our case, polytope or PL ball) to map $\Phi$ from its $r$-fold deleted join, see the book \cite{matousek2003using} for detailed exposition of these techniques.

For the rest of the text we fix  $N=(d+1)(r-1)\geq 1$ and $P$ from the statement of the main theorem.
Faces of polytopes and cells of complexes are assumed to be closed unless specifically said otherwise.
The empty set $\varnothing$ and the polytope itself are considered faces.

In the polytopal case we set
\[
Q=P^{*r}=P_1*\cdots*P_r,
\]
the $r$-fold join of $P$, and 
\[
D=P^{*r}_\Delta,
\]
the $r$-fold deleted join of $P$, that is, the union of $F_1*\dots*F_r\subseteq Q$ over all $r$-tuples of pairwise disjoint faces of $P$. 

In the simplicial-ball case we set
\[
Q=P^{*r}=P_1*\cdots*P_r,
\]
the $r$-fold join of $P$, and 
\[
D=(\partial P)^{*r}_\Delta.
\]

Theorem~\ref{theorem:main} is the standard consequence of a more general Theorem~\ref{theorem:zero}. Before stating it we repeat the following standard construction.

The continuous affine on faces map $\Phi:Q\longrightarrow \R^N$ is defined in the following standard way. Choose vectors $u_1,\ldots,u_r\in\R^{r-1}$ which span
$\R^{r-1}$ and satisfy
\[
u_1+\cdots+u_r=0,
\]
with this relation spanning the space of all linear dependencies among
the $u_i$.  For example, one may take the vertices of a regular
simplex centered at the origin.

Define
\[
\Phi:Q\longrightarrow
\R^{d+1}\otimes\R^{r-1}\cong\R^N
\]
by defining its restriction to the $i$th join factor as
\begin{equation}
\label{eq:Phi}
\Phi_i(x)=(\phi(x),1)\otimes u_i.
\end{equation}

Clearly
\[
\Phi_1(x)+\cdots+\Phi_r(x)
=(\phi(x),1)\otimes(u_1+\cdots+u_r)=0.
\]

\begin{theorem}
\label{theorem:zero}
In the simplicial-ball case, let $\Phi:Q\to\R^N$ be the map defined above in \eqref{eq:Phi} and $r\ge 3$. In the polytopal case, let $\Phi:Q\to\R^N$ be an arbitrary affine map satisfying  $\Phi_1(x)+\cdots+\Phi_r(x)=0$ for all $x\in P$.

Then $0\in\Phi(D)$.
\end{theorem}

The simplicial-ball case of Theorem~\ref{theorem:main} is already known from \cite{hktt2023} for $r$ a prime power and for $r=2$ in particular. So, the following deduction assumes $r\geq 3$ in the simplicial-ball case.

\begin{proof}[Deduction of Theorem~\ref{theorem:main} from Theorem~\ref{theorem:zero}]
In both simplicial-ball and polytopal case let $\Phi:Q\longrightarrow \R^N$ be the map defined above by \eqref{eq:Phi}.

By Theorem~\ref{theorem:zero}, there is a point
\[
z=t_1x_1\oplus\cdots\oplus t_rx_r\in D
\]
such that $\Phi(z)=0$.  Thus
\[
\sum_{i=1}^r t_i(\phi(x_i),1)\otimes u_i=0.
\]

The kernel of the linear map
\[
(\R^{d+1})^r\longrightarrow
\R^{d+1}\otimes\R^{r-1},
\quad
(a_1,\ldots,a_r)\longmapsto\sum_{i=1}^r a_i\otimes u_i,
\]
is the diagonal subspace.  Indeed, this follows coordinate-wise from
the fact that the kernel of
\[
\R^r\longrightarrow\R^{r-1},
\quad
(s_1,\ldots,s_r)\longmapsto\sum_{i=1}^r s_i u_i,
\]
is spanned by $(1,\ldots,1)$.  Consequently,
\[
t_1(\phi(x_1),1)=\cdots=t_r(\phi(x_r),1).
\]
Comparing the last coordinates and using $t_1+\cdots+t_r=1$, we obtain
\[
t_1=\cdots=t_r=\frac1r.
\]
In particular, all the join coefficients are positive, and comparison
of the first $d$ coordinates gives
\[
\phi(x_1)=\cdots=\phi(x_r).
\]

Since $z\in D$ and every $t_i$ is positive, there are pairwise
disjoint nonempty faces $F_i$ of $P$ such that $x_i\in F_i$. Since $r\geq 2$, we have that $F_i\neq P$ for all $i$ and so all $F_i$ belong to $\partial P$. Their
images contain the common point $\phi(x_1)=\cdots=\phi(x_r)$.
\end{proof}

\section{Proof of Theorem~\ref{theorem:zero} modulo lemmas}

All homology and cohomology groups are with coefficients in a fixed field $\F$ which will be omitted from the notation when it is not needed. Vanishing of $\widetilde H_{-1}(\cdot)$ or $\widetilde H^{-1}(\cdot)$ means nonemptiness.

The main technical tool we need is the following standard topological lemma:

\begin{lemma}[Leray--Vietoris--Begle]
\label{lemma:proper-map}
Let $f:A\to B$ be a semialgebraic surjection between compact semialgebraic
spaces.
\begin{enumerate}
\item If every fiber of $f$ is cohomologically acyclic, then
\[
f^*:H^*(B;\F)\longrightarrow H^*(A;\F)
\]
is an isomorphism.
\item Let $n\geq0$.  If every fiber is nonempty and
\[
\widetilde H^j(f^{-1}(b);\F)=0
\quad 0\leq j\leq n-1,
\]
then
\[
f^*:\widetilde H^n(B;\F)
\longrightarrow\widetilde H^n(A;\F)
\]
is injective.
\end{enumerate}
\end{lemma}

\begin{proof}
We use the Leray spectral sequence \cite[Theorem~4.17.1]{godement1958} for the constant sheaf $\F_A$:
\[
E_2^{p,q}=H^p(B;R^qf_*\F_A)
\quad\Longrightarrow\quad
H^{p+q}(A;\F).
\]
Since $f$ is a continuous map from a compact space, the preimages of a base of neighborhoods of any point $b\in B$ produce a base of neighborhoods of the fiber $f^{-1}(b)$ (\emph{proper base change}). Indeed, for any neighborhood $U\supseteq f^{-1}(b)$ its complement $A\setminus U$ and its image $f(A\setminus U)$ are compact and closed, hence $V=B\setminus f(A\setminus U)$ is an open neighborhood of $b$ such that $f^{-1}(V)\subseteq U$.

Hence the stalk of
$R^qf_*\F_A$ at $b\in B$ equals
\[
H^q(f^{-1}(b);\F)
\]
from the continuity of sheaf cohomology \cite[Theorem~4.11.1]{godement1958}. For compact semialgebraic spaces, like $A$, $B$, and the fiber $f^{-1}(b)$, sheaf cohomology with constant coefficients agrees with singular cohomology (this follows from \cite[Theorem~9.2.1]{rag1998}), so we do not distinguish between the kinds of cohomology here.

If every fiber has zero reduced cohomology, then the natural map
\[
\F_B\longrightarrow f_*\F_A
\]
is an isomorphism and $R^qf_*\F_A=0$ for $q>0$.  The spectral sequence
is therefore concentrated in the row $q=0$ and gives
\[
H^*(B;\F)\cong H^*(A;\F).
\]
The resulting isomorphism is the homomorphism $f^*$.

Now suppose that the hypotheses of the second part hold and that
$n\geq1$.  The fibers are connected, so again
\[
f_*\F_A\cong\F_B,
\]
and the fiberwise vanishing and proper base change give
\[
R^qf_*\F_A=0
\quad\text{for any}\quad 1\leq q\leq n-1.
\]
In the Leray spectral sequence, no differential leaves
$E_r^{n,0}$, and every possible differential entering it has source
$E_r^{n-r,r-1}$, which vanishes for $2\leq r\leq n$; for larger $r$
the first index is negative.  Thus
\[
E_\infty^{n,0}=E_2^{n,0}=H^n(B;\F).
\]
The edge homomorphism
\[
f^*:H^n(B;\F)\longrightarrow H^n(A;\F)
\]
is therefore injective.  For $n>0$, this is also the claimed statement
in reduced cohomology.  If $n=0$, injectivity of
$\widetilde H^0(B;\F)\to\widetilde H^0(A;\F)$ follows directly from
the surjectivity of $f$.
\end{proof}

\begin{proof}[Proof of Theorem~\ref{theorem:zero} modulo lemmas]
Assume for contradiction that $0\notin\Phi(D)$.

In the polytopal case we will need $\Phi$ to be non-degenerate, meaning that $\aff\Phi(Q)=\R^N$. We may assume this because of the following perturbation argument.

Define the diagonal map $\Delta : \mathbb R^N\to (\mathbb R^N)^{*r}$ by $\Delta(x) = \frac{1}{r}x\oplus\dots\oplus\frac{1}{r}x$. Choose the origin in the affine span of $Q$ in $\Delta(\mathbb R^N)$ and consider $\Delta(\mathbb R^N)$ as a linear subspace. If $\Phi(D)\not\ni 0$ then we perturb  $\Phi$ slightly as a linear map from $\aff(Q)/\Delta(\mathbb R^N)$. This keeps both $\Phi(D)\not\ni 0$ and the assumption 
\begin{equation}
\label{eq:diagonal-relation}
\Phi(\Delta(P)) = 0\Leftrightarrow\Phi_1(x)+\cdots+\Phi_r(x)=0
\quad\text{for any}\quad x\in P.
\end{equation}
Since $\dim \aff(Q)/\Delta(\mathbb R^N) = (r-1)(N+1) > N$ in our range $r\ge 2$, $d\ge 0$, a generic perturbation surjectively maps $\aff(Q)$ to the whole $\mathbb R^N$.

For any $y\in S^{N-1}$ let the affine (in the simplicial-ball case, affine on faces) functional $\lambda_y:Q\to \R$ be defined by
\[
\lambda_y=\langle y,\Phi(\,\cdot\,)\rangle,
\]
where $\langle\cdot,\cdot\rangle$ is the inner product of the Euclidean space.

Denote
\[
X=
\bigcup_{G\text{ a face of }D}
G\times
\left\{y\in S^{N-1}:
\lambda_y\geq 0\text{ on }G\right\}.
\]

Let
\[
{\rm pr}_D:X\longrightarrow D
\]
be the first projection.  If $x\in D$ and $G_x$ is the minimal face of
$D$ containing $x$, then
\[
{\rm pr}_D^{-1}(x)=
\left\{y\in S^{N-1}:
\lambda_y(z)\geq0
\text{ for every }z\in G_x\right\}.
\]

There is $y_0\in S^{N-1}$ such that 
\[
\lambda_{y_0}(z)>0
\quad\text{for any}\quad z\in G_x.
\]
Indeed, $\Phi(G_x)$ is compact, convex, and does not contain the origin, hence $y_0$ is given by strict separation.



The fiber ${\rm pr}_D^{-1}(x)$ is contractible.  Indeed, for every
$y$ in this fiber and every $t\in[0,1]$, the vector
\[
(1-t)y+ty_0
\]
is nonzero, and its normalization remains in the fiber.  This gives a
contraction of the fiber to $y_0$.

The space $X$ is compact and semialgebraic (defined by inequalities on the finite number of vertices of $Q$), and ${\rm pr}_D$ is a semialgebraic surjection. Lemma~\ref{lemma:proper-map} therefore gives
\[
\widetilde H^{N-1}(X)
\cong\widetilde H^{N-1}(D).
\]

\begin{lemma}[Acyclicity of the deleted join]
\label{lemma:diagonal_acyclicity}
Under the assumption that either $r\ge 3$ or we are considering the polytopal case, $\widetilde H^{N-1}(D)=0.$
\end{lemma}

From the lemma we now get 
\[
\widetilde H^{N-1}(X) = 0.
\]

\begin{lemma}
\label{lemma:Y}
There is a closed semialgebraic set $Y$ with the following properties:
\begin{itemize}
\item[$(I)$] $Y$ is contained in $X$ and the second projection
\[
{\rm pr}_{S^{N-1}}:Y\longrightarrow S^{N-1}
\]
is a surjection.
\item[$(II)$] For any $y\in S^{N-1}$ and $Y_y=({\rm pr}_{S^{N-1}})^{-1}(y)$ we have
\[
\widetilde H^j(Y_y)=0
\quad\text{for any}\quad j\leq N-2.
\]
\end{itemize}
\end{lemma}

Lemma~\ref{lemma:proper-map}, with $n=N-1$, implies that
\[
{\rm pr}_{S^{N-1}}^*:
\widetilde H^{N-1}(S^{N-1})
\longrightarrow\widetilde H^{N-1}(Y)
\]
is injective. Since
$\widetilde H^{N-1}(S^{N-1})\cong\F$, this map is nonzero.

Let $i:Y\hookrightarrow X$ be the inclusion, and let
\[
{\rm pr}_{S^{N-1}}:X\longrightarrow S^{N-1}
\]
also denote the second projection. The second projection on $Y$ is
the composite
\[
Y\xrightarrow{i}X\xrightarrow{{\rm pr}_{S^{N-1}}}S^{N-1}.
\]
But $\widetilde H^{N-1}(X)=0$, so the homomorphism induced by this
composite is zero.  This contradicts the injectivity proved above.
Therefore $0\in\Phi(D)$.

\end{proof}

\section{Proof of the simplicial-ball case of Lemma~\ref{lemma:Y}}

\begin{proof}[Proof of the simplicial-ball case of Lemma~\ref{lemma:Y}]
Let $L\subset \partial P$ be an $(N-1)$-dimensional face, and define
\[
Y=
\bigcup_{G\text{ a face of }L^{*r}_\Delta}
G\times
\left\{y\in S^{N-1}:
\lambda_y\geq0\text{ on }G\right\}.
\]
Every face of $L^{*r}_\Delta$ is a face of $D$, so the definition of
$X$ immediately gives $Y\subset X$.  The set $Y$ is compact and
semialgebraic (speaking of semialgebraic here we may assume $P$ affinely embedded into some $\mathbb R^M$ of high dimension), and the second projection
\[
{\rm pr}_{S^{N-1}}:Y\longrightarrow S^{N-1}
\]
is algebraic.

Let $v_1,\ldots,v_N$ be the vertices of $L$.  For $y\in S^{N-1}$ and a
vertex $v$ of $L$, let
\[
W_v(y)=
\left\{v^{(i)}\in P_i:
\lambda_y(v^{(i)})\geq0,
\quad 1\leq i\leq r\right\},
\]
where $v^{(i)}$ denotes the copy of $v$ in $P_i$.  The diagonal
relation gives
\[
\sum_{i=1}^r\lambda_y(v^{(i)})
=
\left\langle y,\sum_{i=1}^r\Phi_i(v)\right\rangle
=0.
\]
Hence $W_v(y)$ is nonempty for every vertex $v$ of $L$.

Since $\lambda_y$ is affine on every face of $L^{*r}_\Delta$, its
nonnegativity on a face is equivalent to its nonnegativity on all the
vertices of that face.  It follows that the fiber over $y$ is
\[
Y_y=
W_{v_1}(y)*\cdots*W_{v_N}(y).
\]
In particular, every fiber is nonempty, so the second projection is
surjective.  A join of $N$ nonempty spaces is $(N-2)$-connected when
each factor is a nonempty discrete space.  Therefore
\[
\widetilde H^j(Y_y)=0
\quad\text{for every}\quad j\leq N-2.
\]
Thus $Y$ satisfies both parts of Lemma~\ref{lemma:Y}.
\end{proof}

\section{The polytope $R_y$, polytopal $Y$, and the proof of part (I) of Lemma~\ref{lemma:Y}}
\label{section:Ry}

From now on we work in the case of Theorem~\ref{theorem:zero} where $P$ is a polytope and both $\phi$ and $\Phi$ are affine maps.

Consider the projection to the diagonal $\pi:Q\longrightarrow P$,
\[
\pi(t_1x_1\oplus\cdots\oplus t_rx_r)
=t_1x_1+\cdots+t_rx_r.
\]
Introduce the function $h_y: P\to \mathbb R$,
\[
h_y(x)=
\max\{\lambda_y(w):w\in Q,\ \pi(w)=x\}.
\]
Consider the following affine image of $Q$:
\[
R_y=(\pi,\lambda_y)(Q) \subset P\times\mathbb R,
\]
and denote the top part of $R_y$ (seen from the positive direction of the last coordinate) by
\[
\Gamma_y=\{(x,h_y(x)):x\in P\}.
\]
It follows that $h_y$ is a continuous concave function on $P$.

Recall that by the argument in the proof of Theorem~\ref{theorem:zero} we may assume $\Phi$ non-degenerate, meaning that $\aff\Phi(Q)=\R^N$.

\begin{claim}\label{claim:R_y-full-dimensional}
The dimension of the polytope $R_y$ is $N+1$.
\end{claim}
\begin{proof}
The projection of $R_y$ onto $P$ is surjective, so its dimension is at
least $N$.  If its dimension were $N$, projection would map
$\aff R_y$ affinely isomorphically onto
$\aff P$.  Hence $R_y$ would be contained in the graph of
an affine function $\beta:P\to\R$, and
\[
\lambda_y=\beta\circ\pi
\quad\text{on}\quad Q.
\]
The restriction of this equality to each factor gives
$\langle y,\Phi_i(\,\cdot\,)\rangle = \lambda_y|_{P_i}=\beta$.
From \eqref{eq:diagonal-relation} we now get $r\beta=0$, so $\lambda_y=0$ on $Q$.
This would mean that $\Phi(Q)$ lies in the orthogonal complement to $y$, which contradicts  the non-degeneracy of $\Phi$.
\end{proof}

We are ready to define $Y$ which we use in the proof of the polytopal case of Lemma~\ref{lemma:Y}. Put:
\[
Y=
\left\{(z,y)\in D\times S^{N-1}:
\lambda_y(z)=h_y(\pi(z))
\right\}.
\]

The fiber $Y_y$ over $y\in S^{N-1}$ under the second projection $Y\longrightarrow S^{N-1}$ satisfies
\[
Y_y=\{z\in D : \lambda_y(z)=h_y(\pi(z))\},
\]
and
\begin{equation}\label{eq:Y_y-formula}
Y_y=D\cap (\pi,\lambda_y)^{-1}(\Gamma_y).
\end{equation}

For every nonempty face $F$ of $Q$, let $S_F$ be the set of all
$y\in S^{N-1}$ for which there is an affine function
$\gamma:P\to\R$ such that
\begin{equation}\label{eq:SF-inequality}
\lambda_y(w)\leq\gamma(\pi(w))\quad \text{for any}\quad w\in Q \text{\quad and \quad } \lambda_y(z)=\gamma(\pi(z)) \quad \text{for any}\quad z\in F.
\end{equation}

\begin{claim}\label{claim:Y-face-description}
For every $y\in S^{N-1}$, we have equality
\[
Y_y=
\bigcup_{F\text{ a nonempty face of }Q,\ y\in S_F}
(D\cap F).
\]
\end{claim}

\begin{proof}
Suppose first that $y\in S_F$ and let us prove that $D\cap F\subset Y_y$. Fix $z\in D\cap F$ and let $\gamma$ be a functional satisfying \eqref{eq:SF-inequality}. Let $x = \pi(z)$. For any other $w\in Q$ with $\pi(w)=x$ we have 
\[
\lambda_y(w)\leq\gamma(\pi(w))=\gamma(x)=\gamma(\pi(z))=\lambda_y(z)
\]
by \eqref{eq:SF-inequality}. Hence, $\lambda_y(z) = \max\{\lambda_y(w):w\in Q,\ \pi(w)=x\} = h_y(x)$ and $z\in Y_y$ by the definition of $Y_y$.

Conversely, let $z\in Y_y$ and put
\[
p=(\pi(z),\lambda_y(z))\in\Gamma_y.
\]
By the definition of $R_y$, no point ``above'' $p$ belongs to $R_y$, i.e., there is no $s>0$ such that $p+se\in R_y$, where $e=(0,1)$ is the ``up'' vector. Take one of the finite number of linear inequalities $ax+bt+c\le 0$ defining the polytope $R_y$ which is violated by $p+se$ when $s$ tends to $0$, it must have $b>0$ and be an equality at $p$. Hence, there is an affine functional
\[
(x,t)\longmapsto ax+bt+c
\]
which is non-positive on $R_y$, zero at $p$, and positive on $p+se$ for any $s>0$. In particular,
$b>0$, and this functional is maximized on $R_y$ at $p$.

Define an affine functional $\gamma:P\to\R$ by
\[
\gamma(x)=-\frac{ax+c}{b}
\]
and $F\subset Q$ by
\[
F = \{w\in Q: \lambda_y(w) = \gamma({\pi(w)})\}.
\]

It remains to check that $\gamma$ satisfies \eqref{eq:SF-inequality} and that $F\subset Q$ is a face with $z\in F$.

For every $w\in Q$ we have $(\pi(w),\lambda_y(w))\in R_y$ by definition.
Since the supporting functional $(x,t)\longmapsto ax+bt+c$ is nonpositive on
$R_y$ and zero at $p$,
\[
a\pi(w)+b\lambda_y(w)+c\leq 0
\]
with the equality achieved if (but not only if) $w=z$.
Because $b>0$,
\[
\lambda_y(w)
\leq-\frac{a\pi(w)+c}{b}
=\gamma(\pi(w))
\]
and \eqref{eq:SF-inequality} is satisfied.

The affine function
\[
w\longmapsto
\lambda_y(w)-\gamma(\pi(w))
\]
is nonpositive on $Q$, $F$ is its maximizing set and it is maximized at $z$.  Therefore $F\ni z$ is a face of $Q$.

\end{proof}

\begin{claim}\label{claim:Yy-nonempty}
Every fiber $Y_y$ is nonempty.
\end{claim}
\begin{proof}
The affine function $\lambda_y$ attains its maximum on $Q$ on a nonempty face of $Q$.  Choose a vertex $v$ of this face.  Every vertex of $Q$ belongs to the deleted join $D$.  Moreover, if $w\in Q$ satisfies $\pi(w)=\pi(v)$, then
\[
\lambda_y(w)\leq\lambda_y(v),
\]
because $v$ is a global maximizer of $\lambda_y$ on $Q$.  Therefore
\[
\lambda_y(v)=h_y(\pi(v)),
\]
and hence $v\in Y_y$.
\end{proof}

\begin{claim}\label{claim:Y-closed}
The set $Y$ is closed in $D\times S^{N-1}$.
\end{claim}

\begin{proof}
Identify the space of affine functionals on $\R^N$ with $\R^{N+1}$.
Let $F$ be a face of $Q$.
Define $C_F\subset \R^N\times \R^{N+1}$ as the set of pairs $(y,\gamma)$ ($\gamma$ considered as an affine function on $\mathbb R^N$) such that
\begin{align*}
\langle y,\Phi(w)\rangle = \lambda_y(w)
&\leq \gamma(\pi(w))
&&\text{for every vertex }w\in Q,\\
\langle y,\Phi(z)\rangle = \lambda_y(z)
&=\gamma(\pi(z))
&&\text{for every vertex }z\in F.
\end{align*}

The set $C_F$ is precisely the set of pairs $(y,\gamma)$ such that $\gamma$ satisfies \eqref{eq:SF-inequality} with the given $y$ on the vertices of $Q$ and $F$, and hence on the whole $Q$ and $F$. So, $S_F$ is the intersection of $S^{N-1}$ with the projection of $C_F$ to the first coordinate. Both $C_F$ and its projection are polyhedra, hence $S_F$ is closed as the intersection of $S^{N-1}$ with a (possibly unbounded) polyhedron.

So, every term in the finite union representing $Y$ in Claim~\ref{claim:Y-face-description} is closed. Therefore, $Y$ itself is closed.
\end{proof}

\begin{claim}\label{claim:Y-semialgebraic}
The set $Y$ is semialgebraic.
\end{claim}

\begin{proof}
For every face $F$ of $Q$, the set $D\cap F$ is semialgebraic.  The set
$S_F$ is the intersection of the sphere with a polyhedron and is
therefore semialgebraic. By Claim~\ref{claim:Y-face-description},
$Y$ is a finite union of products of such sets.  Hence $Y$ is
semialgebraic.
\end{proof}

\begin{claim}\label{claim:Y-in-X}
We have $Y\subset X$.
\end{claim}

\begin{proof}
By Claim~\ref{claim:Y-face-description}, it is enough to check that if $F\subset Q$ is a face and $y\in S_F$ then $\lambda_y(F)\geq 0$. 

Let $z$ be a vertex of $F$, it corresponds to a vertex $\pi(z)$ of $P$. By~\eqref{eq:diagonal-relation}, there exists a vertex $w$ of $Q$ (possibly $w=z$) which corresponds to the same vertex $\pi(w)=\pi(z)$ of $P$ and such that $\langle y,\Phi(w)\rangle=\lambda_y(w)\geq 0$.

By the definition of $S_F$~\eqref{eq:SF-inequality}, we have 
\[
\lambda_y(z)=\gamma(\pi(z))=\gamma(\pi(w))\geq \lambda_y(w) \geq 0.
\]

So, $\lambda_y$ is non-negative on all the vertices of $F$ and hence on the whole $F$.
\end{proof}

\begin{proof}[Proof of part (I) of polytopal Lemma~\ref{lemma:Y}]
We have $Y\subset X$ by Claim~\ref{claim:Y-in-X}. The set $Y$ is semialgebraic by Claim~\ref{claim:Y-semialgebraic}. The projection ${\rm pr}_{S^{N-1}}:Y\to S^{N-1}$ is surjective by Claim~\ref{claim:Yy-nonempty}. Finally, it is closed by Claim~\ref{claim:Y-closed}.
\end{proof}

\section{Bad-vertex triangulation lemmas and the proof of Lemma~\ref{lemma:diagonal_acyclicity}}


In the polytopal case for a nonempty face $F=Q_1*\cdots*Q_r$ of $Q$ we put
\[
m(F)=\dim\aff (\pi(Q_1)\cup\cdots\cup \pi(Q_r)).
\]
In the simplicial-ball case for a nonempty face $F=Q_1*\cdots*Q_r$ of $(\partial P)^{*r}$ we denote by $m(F)$ the number of distinct vertices of $P$ in $Q_1\cup\cdots\cup Q_r\subset\partial P$ minus $1$.

The simplicial-ball case of the following lemma and its proof is obtained by replacing $Q$ with $(\partial P)^{*r}$ and $D$ with $(\partial P)^{*r}_\Delta$.

\begin{lemma}[Bad-vertex triangulation]\label{lemma:bad-vertex-triangulation}
There is a triangulation $\mathcal T$ of $Q$, together with a
designation of its vertices as good or bad, such that the following
properties hold.
\begin{enumerate}
\item The subcomplex induced by the good vertices triangulates $D$.
\item For every nonempty face $F=Q_1*\cdots*Q_r$ of $Q$, every simplex
of $\mathcal T|_F$ which is full-dimensional in $F$ contains at least
$m(F)+1$ good vertices.
\end{enumerate}
\end{lemma}

\begin{proof}
We triangulate the skeleta of $Q=P^{*r}$ by induction over their
dimension.  Call an edge of $Q$ bad if its endpoints are two copies of
the same vertex of $P$.  A face of $Q$ belongs to the deleted join $D$ if
and only if it contains no bad edge.

Choose an arbitrary order on the vertices of $Q$ and on its bad edges.
Suppose that the $(k-1)$-skeleton has already been triangulated, and
let $F$ be a $k$-dimensional face.  If $F$ belongs to the deleted join $D$,
triangulate it by coning the already triangulated boundary from the
earliest available vertex of $F$.  If $F$ does not belong to the
deleted join, triangulate it by coning from the midpoint of the
earliest available bad edge of $F$.  In either case, only the boundary
faces in $\partial F$ which do not contain the chosen apex $v$ are used in the coning.
As for the faces in $\partial F$ which do contain $v$, 
our chosen fixed order on the vertices of $Q$ and on its bad edges guarantees that they have already been triangulated at the earlier steps using the same apex $v$. This ensures that the triangulation constructed at the $k$th step is compatible with the previous steps. 

Declare the midpoints of bad edges to be bad vertices of the resulting
triangulation, and declare all other vertices to be good.  Every face
belonging to the deleted join is triangulated using only good vertices.
Conversely, every full-dimensional simplex subdividing a face which
does not belong to the deleted join contains the bad apex used to
triangulate that face.  Induction over the dimension of the minimal
face containing a simplex therefore shows that a simplex contains only
good vertices if and only if it lies in $D$.  Thus the subcomplex
induced by the good vertices triangulates $D$.

It remains to prove the second property.  Consider a
full-dimensional simplex $S$ of the triangulation of a face
\[
F=Q_1*\cdots*Q_r.
\]
Let $k=\dim F$ and $v_k$ be the vertex of $S_k=S$ which was used as an apex when triangulating $F_k=F$.
Let $S_{k-1}$ be the facet of $S_k$ opposite $v_k$ then $S_{k-1}$ is a full-dimensional simplex of some facet $F_{k-1}$ of $F_k$. Iterating this reasoning we get a flag 
\[
F_0\subset F_1\subset\cdots\subset F_k=F,
\]
where $\dim F_i=i$ and the $i$th vertex $v_i$ of $S$ is the apex used
when triangulating $F_i$.

We have that
\[
m(F_{i-1})\leq m(F_i)\leq m(F_{i-1})+1,
\quad
m(F_0)=0,
\]
and that $v_0$ is good.

We claim that if $v_i$ is bad, then
\[
m(F_i)=m(F_{i-1}).
\]

Indeed, write
\[
F_\ell=Q_{1,\ell}*\cdots*Q_{r,\ell},
\] 
for any $\ell$. After relabeling the join factors, we may assume that
$Q_{1,i}$ is a coface of $Q_{1,i-1}$ and that
$Q_{j,i}=Q_{j,i-1}$ for $j\geq2$.  

Since $v_i$ is bad, the bad edge
chosen in $F_i$ is not contained in $F_{i-1}$.  Hence it has one
endpoint which is a vertex of
$Q_{1,i}\setminus Q_{1,i-1}$, while a copy of the same vertex of $P$
already belongs to some $Q_{j,i}=Q_{j,i-1}$ with $j\geq2$.  Thus the
new point added to the affine hull of $Q_{1,i-1}$ already lies in
\[
\aff 
\bigl(\pi(Q_{1,i-1})\cup\cdots\cup\pi(Q_{r,i-1})\bigr),
\]
which proves the claim.

Consequently, the integer $m(F_i)$ can increase only when $v_i$ is
good.  It starts at $0$ and ends at $m(F)$, and $v_0$ is also good.
Therefore $S$ has at least $m(F)+1$ good vertices.
\end{proof}

The polytopal case of Lemma~\ref{lemma:diagonal_acyclicity} follows from the special case $F=Q$ of the following lemma and the fact that the coefficients take values in a field.

\begin{lemma}[Local acyclicity]\label{lemma:local-acyclicity}
For every nonempty face $F$ of $Q$,
\[
\widetilde H_j(D\cap F)=0
\quad\text{for any}\quad j\leq m(F)-1.
\]
\end{lemma}

\begin{proof}
Put $k=\dim F$ and let $L=\mathcal T|_F$.  Let $K\subset L$ and
$B\subset L$ be the subcomplexes induced by the good and bad vertices,
respectively.  By Lemma~\ref{lemma:bad-vertex-triangulation},
\[
D\cap F = |K|.
\]

Every $k$-simplex of $L$ contains at least $m(F)+1$ good vertices and
therefore at most $k-m(F)$ bad vertices.  Every simplex of $B$ is
contained in a $k$-simplex of $L$, so
\[
\dim B\leq k-m(F)-1.
\]

The set of vertices of $L$ is partitioned into the set of vertices of $K$ and of $B$ and the subcomplexes $K$ and $B$ are induced. Hence, the complement $|L|\setminus|B|$ deformation retracts onto $|K|$. 

Fix $0\leq j\leq m(F)-1$. Since $L$ is a ball, any $j$-cycle in $K$ can be filled by a $(j+1)$-dimensional PL chain in $L$. Move the chain into general position with $|B|$ while keeping it fixed on its boundary -- the cycle. After that the intersection of the chain with $|B|$ has dimension at most $(j+1)+\dim B-\dim L = (j+1)+(k-m(F)-1)-k=j-m(F)<0$, meaning that the chain avoids $|B|$ and lies in $|L|\setminus|B|$. Retracting the chain back to $|K|$ we get that the cycle was trivial, proving
\[
\widetilde H_j(|L|\setminus|B|)=0
\quad\text{for any}\quad 0\leq j\leq m(F)-1.
\]
The complex $K$ is nonempty because every maximal simplex of $L$
contains a good vertex, so the same conclusion also holds in degree
$j=-1$.
\end{proof}

To prove the simplicial-ball case of Lemma~\ref{lemma:diagonal_acyclicity} a bit of additional work is required.

\begin{proof}[Proof of the simplicial-ball case of Lemma~\ref{lemma:diagonal_acyclicity}]
Let $K$ and $B$ be the subcomplexes of $\mathcal T$ induced by the good and bad vertices,
respectively. By Lemma~\ref{lemma:bad-vertex-triangulation}, $K$ triangulates $(\partial P)^{*r}_\Delta$.

Since $(\partial P)^{*r}\cong S^{rN-1}$, Alexander duality together with the fact that the complement to $|K|$ retracts to $|B|$ gives
\[
\widetilde H^{N-1}(K)\cong\widetilde H_b((\partial P)^{*r}\setminus |K|)\cong\widetilde H_b(B),
\quad\text{where}\quad b=(r-1)N-1.
\]

Let $\tau$ be a $b$-simplex of $B$. It is a face of some maximal simplex $T$ of $\mathcal T$, which appeared in the construction when we were triangulating a maximal face $\sigma_1*\ldots *\sigma_r$ with each $\sigma_i$ a facet of $\partial P$.

We have that $m(\sigma_1*\ldots *\sigma_r)\geq N-1$ since each individual $\sigma_i$ already has $N$ distinct vertices. By Lemma~\ref{lemma:bad-vertex-triangulation}, $T$ has at least $m(\sigma_1*\ldots *\sigma_r)+1\geq N$ good vertices. On the other hand, $T$ has at most $N$ good vertices, since all the vertices in $\tau$ are bad. We get $m(\sigma_1*\ldots *\sigma_r)=N-1$, meaning that $\sigma_1,\ldots ,\sigma_r$ collectively have only $N$ distinct vertices and implying $\sigma_1=\ldots=\sigma_r=\sigma$.

Let $v$ be the bad apex used when triangulating $\sigma^{*r}$. Then $v$ is a vertex of $\tau$. Let $\tau_v\subset\tau$ be the opposite facet.
We now want to show that there is no other $b$-simplex $\tau'\neq\tau$ of $B$ having $\tau_v$ as its facet. Analogously to $\tau$, $\tau'$ must have appeared when we were triangulating $\sigma'^{*r}$ for some facet $\sigma'$ of $\partial P$.

Assume first that $\sigma'=\sigma$. Since $\tau'\neq \tau$ we have that $v\not\in \tau'$. Let $T'$ be a maximal simplex in the triangulation of $\sigma^{*r}$ which contains $\tau'$. By Lemma~\ref{lemma:bad-vertex-triangulation}, $T'$ contains at least $m(\sigma^{*r})+1=N$ good vertices; it also contains $b+2$ bad vertices, $v$ plus $b+1$ vertices from $\tau'$. Summing up we get at least $N+b+2=rN+1$ vertices, while the dimension of $T'$ is only $\dim \sigma^{*r}=rN-1$, contradiction.

Assume now that $\sigma'\neq\sigma$. Let $a$ be the number of vertices in $\sigma'\cap\sigma$, we have $a\leq N-1$ because $\sigma'$ and $\sigma$ are different. By Lemma~\ref{lemma:bad-vertex-triangulation}, every maximal simplex of  $\mathcal T$ within $(\sigma'\cap\sigma)^{*r}$ contains at least $m((\sigma'\cap\sigma)^{*r})+1=a$ good vertices and hence at most $\dim (\sigma'\cap\sigma)^{*r}+1-a=a(r-1)$ bad.

On the other hand, $\dim \tau_v=b-1$ so it has $(r-1)N-1$ bad vertices. 
Since we assume $r\geq 3$, we have that $a(r-1)\leq (N-1)(r-1) < (r-1)N-1$, so $\tau_v$ has more bad vertices than any maximal simplex of $\mathcal T$ within $(\sigma'\cap\sigma)^{*r}$, which contradicts $\tau_v\subset (\sigma'\cap\sigma)^{*r}$.

So, we have proved that any bad $b$-simplex $\tau$ has a free facet $\tau_v$.
Hence, no non-trivial $b$-cycles can exist in $B$. From that we conclude that
$\widetilde H_b(B) = 0$ and $\widetilde H^{N-1}(D)=\widetilde H^{N-1}(K)=0$.
\end{proof}

\section{Shellability and acyclicity lemmas}

To study the acyclicity further we need two more lemmas.
The second one is the classical polytope shellability result.

\begin{lemma}[Acyclic gluing over a shellable ball]\label{lemma:acyclic-gluing}
Let $\mathcal S$ be a shellable decomposition into convex polytopal cells of an $n$-dimensional PL ball. Let $K$ be a finite simplicial complex. Suppose that to every cell $E$ of $\mathcal S$ we assign a subcomplex $A_E\subset K$ such that
$A_\varnothing=\varnothing$ and
\[
A_E\cap A_{E'}=A_{E\cap E'}
\]
for all cells $E,E'$.  Suppose also that, for every nonempty cell $E$,
\[
\widetilde H_j(A_E)=0
\quad\text{for any}\quad j\leq\dim E-1.
\]
Then
\[
\widetilde H_j
\left(\bigcup_{C\text{ maximal in }\mathcal S}A_C\right)=0
\quad\text{for any}\quad j\leq n-1.
\]
\end{lemma}

\begin{proof}
We argue by induction on $n$ and, for fixed $n$, on the number of
maximal cells of $\mathcal S$.  For $n=0$, the ball $\mathcal S$ has
one maximal cell $C$, and the assertion says only that $A_C$ is
nonempty, which is one of the assumptions.

Let $n\geq1$, and let
\[
C_1,\ldots,C_s
\]
be a shelling of the maximal cells of $\mathcal S$.  Put
\[
B_k=A_{C_1}\cup\cdots\cup A_{C_k}.
\]
For $k=1$, the required vanishing for $B_1=A_{C_1}$ is an assumption.
Suppose that it has already been proved for $B_{k-1}$. By the definition of ``shellable'' \cite[Definition~8.1]{ziegler1995} the intersection
\[
\mathcal L_k=C_k\cap(C_1\cup\cdots\cup C_{k-1})
\]
inherits the cellular decomposition whose $(n-1)$-cells form the proper initial segment of a shelling of $\partial C_k$. That is, $\mathcal L_k$ is itself a shellable decomposition of a PL ball into convex polytopal cells, but with smaller dimension.

Using the intersection property of the $A_E$, we obtain
\[
A_{C_k}\cap B_{k-1}
=\bigcup_{i<k}A_{C_k\cap C_i}
=\bigcup_{E\text{ maximal in }\mathcal L_k}A_E.
\]
The assignment $E\mapsto A_E$ restricted to the cells of
$\mathcal L_k$ satisfies the assumptions of the lemma.  Induction on
the dimension therefore gives
\[
\widetilde H_j(A_{C_k}\cap B_{k-1})=0
\quad\text{for any}\quad j\leq n-2.
\]

By induction on $k$ and by the assumption on the maximal cell $C_k$,
\[
\widetilde H_j(B_{k-1})=
\widetilde H_j(A_{C_k})=0
\quad\text{for any}\quad j\leq n-1.
\]
The reduced Mayer--Vietoris sequence for
$B_k=B_{k-1}\cup A_{C_k}$ now implies
\[
\widetilde H_j(B_k)=0
\quad j\leq n-1.
\]
Taking $k=s$ proves the assertion.
\end{proof}

The following lemma is due to Bruggesser and Mani \cite[Proposition 2]{bruggesser1971shellable}, see also \cite[Theorem~8.12]{ziegler1995} and \cite[Claim (3) in the proof of Theorem~15.4]{gruber2007}.

\begin{lemma}[Line shelling]\label{lemma:line-shelling}
Let $R$ be a $d$-dimensional polytope and let $v$ be a generic point
outside $R$. There is a shelling of the facets of $R$ in which the
facets visible from $v$ form an initial segment. Consequently, the
closed face complex generated by the facets visible from $v$ is a
shellable $(d-1)$-dimensional cellular ball.
\end{lemma}

\section{Acyclicity of $Y_y$ and  the proof of the polytopal part (II) of Lemma~\ref{lemma:Y}}

\begin{proof}[Proof of the polytopal part (II) of Lemma~\ref{lemma:Y}]

By Claim~\ref{claim:R_y-full-dimensional} we know that $R_y$ is an $(N+1)$-dimensional polytope. Write every facet inequality of $R_y$ in the form
\[
ax+bt\leq c
\quad\text{for any}\quad (x,t)\in R_y.
\]
Call the facet top, bottom, or vertical according to $b>0$, $b<0$, or
$b=0$.  Let $\mathcal U_y$ be the closed face complex generated by the
top facets.  Its support is exactly $\Gamma_y$.  

Indeed, let us first show that every top facet lies in $\Gamma_y$. Let $(x,t)$ be a point of a top facet satisfying the equality $ax+bt = c$ in the facet's support inequality $ax+bt \leq c$. Then, because $b>0$, every point $(x,t+s)$ with $s>0$ violates the inequality and hence is outside of $R_y$. So, $t=h_y(x)$ and $(x,t)\in\Gamma_y$ by the definitions of $h_y$ and $\Gamma_y$, respectively.

Now we show that every point $p\in\Gamma_y$ lies in a top facet. If all
facets containing $p$ had last coefficient at most zero, then a
sufficiently small movement of $p$ in the direction $(0,1)$ would
still satisfy every facet inequality of $R_y$, which is impossible by the definition $\Gamma_y$.

So, $|\mathcal U_y|=\Gamma_y$. Because $\Gamma_y$ is the graph of a piecewise affine function $h_y$ on $P$, we get that the projections of the faces of $\mathcal U_y$ to $P$ form an
$N$-dimensional polytopal subdivision $\mathcal S_y$ of $P$.  Its
maximal cells are the projections of the top facets.

For every cell $E\in\mathcal S_y$, let $\widetilde F_E$ be the unique
face of $\mathcal U_y$ projecting onto $E$. With the convention $\widetilde F_\varnothing=\varnothing$ we have that
\[
\widetilde F_{E\cap E'}=\widetilde F_E\cap \widetilde F_{E'}
\]
for any two cells $E,E'\in\mathcal S_y$ by the construction of the subdivision $\mathcal S_y$.

Put
\[
F_E=(\pi,\lambda_y)^{-1}(\widetilde F_E) \subset Q.
\]
The set $F_E$ is a face of $Q$ because $(\pi,\lambda_y)(Q)=R_y$ and $(\pi,\lambda_y)$ is affine, hence the preimage of a face $\widetilde F_E$ of $R_y$ is a face of $Q$.
We have 
\[
\pi(F_E)=E.
\]

Define
\[
A_E=D\cap F_E.
\]
From $\widetilde F_{E\cap E'}=\widetilde F_E\cap \widetilde F_{E'}$ and the definition of $F_E$ we get $F_E\cap F_{E'}=F_{E\cap E'}$ and hence
\[
A_E\cap A_{E'}=A_{E\cap E'}.
\]

With respect to the triangulation $\mathcal T$ from
Lemma~\ref{lemma:bad-vertex-triangulation}, every $A_E$ is a subcomplex
of the fixed simplicial complex $\mathcal T|_D$.
If $F_E=Q_1*\cdots*Q_r$, then
\[
E=\pi(F_E)=\conv(\pi(Q_1)\cup\cdots\cup \pi(Q_r)),
\]
and consequently
\[
m(F_E)=\dim E.
\]
Lemma~\ref{lemma:local-acyclicity} gives
\[
\widetilde H_j(A_E)=0
\quad \text{for any}\quad j\leq\dim E-1.
\]

Recall that by \eqref{eq:Y_y-formula} we have that $Y_y=D\cap (\pi,\lambda_y)^{-1}(\Gamma_y)$. Any top facet of $\mathcal U_y$ has the form $\widetilde F_C$, where $C$ is maximal in $\mathcal S_y$. So, $D\cap (\pi,\lambda_y)^{-1}(\widetilde F_C)=D\cap F_C=A_C$. Since the union of the top facets of $\mathcal U_y$ is exactly $\Gamma_y$ we get that
\[
Y_y=
\bigcup_{C\text{ maximal in }\mathcal S_y}A_C.
\]

It remains to verify that $\mathcal S_y$ is shellable.  Choose a generic point 
$x_0\in\operatorname{int}P$.  For sufficiently large and generic $T$, the facets
of $R_y$ visible from the generic point $(x_0,T)$ are exactly the top facets.  Indeed,
the inequalities of the top facets are violated for large $T$, those
of the bottom facets are satisfied, and the inequalities of the
vertical facets are strict because $x_0$ lies in the interior of $P$.
Lemma~\ref{lemma:line-shelling} shows
that $\mathcal U_y$ is a shellable $N$-ball.  Projection identifies
$\mathcal U_y$ with $\mathcal S_y$, so $\mathcal S_y$ is also a
shellable $N$-ball.

Lemma~\ref{lemma:acyclic-gluing} now gives
\[
\widetilde H_j(Y_y)=0
\quad\text{for any}\quad j\leq N-1.
\]
Over a field, the same vanishing holds in reduced cohomology:
\[
\widetilde H^j(Y_y)=0
\quad\text{for any}\quad j\leq N-1.
\]

\end{proof}

\bibliography{../Bib/karasev}

\bibliographystyle{abbrv}

\end{document}